\documentclass [12pt] {article}
\title {The Jones and Alexander polynomials for singular links}
\author{Thomas Fiedler}
\usepackage{amssymb,amsfonts,epsfig}
\begin{document}
\newtheorem{proposition}{Proposition}
\newtheorem{theorem}{Theorem}
\newtheorem{lemma}{Lemma}
\newtheorem{corollary}{Corollary}
\newtheorem{example}{Example}
\newtheorem{remark}{Remark}
\newtheorem{definition}{Definition}
\newtheorem{question}{Question}
\newtheorem{conjecture}{Conjecture}
\maketitle
\begin{abstract}
We extend the Kauffman state models of the Jones and Alexander polynomials of classical links to state models 
of their two-variable extensions in the case of singular links. Moreover, we extend both of them to polynomials with
d+1 variables for long singular knots with exactly d double points. These extensions can detect non-invertibility of
long singular knots.

\footnote{
classical and singular links, knot polynomials , state models}
\end{abstract}

\section{Introduction and results}

{\em Singular links} are oriented smooth links in 3-space but which can have a finite number of rigid ordinary double points.
(compare e.g. \cite{K3}).
The HOMFLY-PT and Kauffman polynomials were extended to 3-variable polynomials for singular links by Kauffman and 
Vogel \cite{KV}. Recently, the extended HOMFLY-PT polynomial was recovered by the construction of traces on singular
Hecke algebras \cite{PR}).

Singular links became popular in connection with the construction of finite type knot invariants (see \cite{V} and 
also \cite{BN}). Recently, there is a new interest in singular links in connection with categorifications.
The Khovanov homology was extended to a homology for singular links in  \cite{S}.
The Alexander polynomials of a cube of resolutions (in Vassiliev's sense) of a singular 
knot were categorified in  \cite{A}. Moreover, a 1-variable extension of the Alexander polynomial for singular links was 
categorified in \cite{OSS}. The generalized cube of resolutions (containing Vassilievs resolutions as well as those 
smoothings at double points which preserve the orientation) was categorified in \cite{OS}.

In this short note we make the simple observation  that  one can easily extend the state models of the Jones and 
Alexander polynomials to state models of their two-variable extensions, called respectively  
$V^s_L(A,B) \in \mathbb{Z}[A ,A^{-1},B]$  and  $\Delta^s_L(A,B) \in \mathbb{Z}[A ,A^{-1},B]$,
in the case of singular links. The extension of the Alexander polynomial in \cite{OSS} coincides with $\Delta^s_L(A,0)$.

The categorification of  $V^s_L(A,B)$ (as a 2-variable polynomial) does not lead to anything new but just several copies of the usual Khovanov
homology \cite{A2}.

\begin{question}
Can $\Delta^s_L(A,B)$ be categorified with a triple graded homology ?
\end{question} 

\begin{question}
Can the state sum from \cite{F} be extended for singular classical links in order to give a 3-variable extension of the
polynomial $W_L$ ?
\end{question}

We assume that the reader is familar with Kauffman's state model for the Jones polynomial \cite{K2} and with his state 
model for the Alexander polynomial \cite{K1}. Moreover, we assume that the link $L$ has a distinguished component, 
which we represent as a long knot. We can then choose the star-regions in a canonical way in the Kauffman state
model of the Alexander polynomial (compare \cite{F}).

Let a singular link be represented by a (singular) link diagram. Two such link diagrams represent the same singular 
link if and only if they can be related by a finite sequence of ordinary Reidemeister moves (see e.g. \cite{BZ}) together 
with the moves
$SII$ shown in Fig. 1 and $SIII$ shown in Fig. 2 (compare \cite{K3}). In particular, the number of double points is an invariant.

\begin{figure}
\centering 
\psfig{file=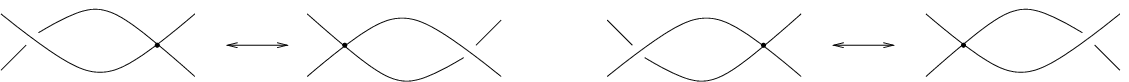}
\caption{}
\end{figure}

\begin{figure}
\centering 
\psfig{file=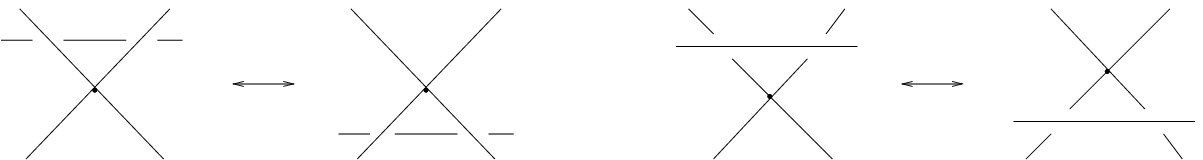}
\caption{}
\end{figure}

Let $p$ be a double point in an oriented singular link diagram $D$.

\begin{definition}
The {\em smoothings} at $p$ are defined in Fig. 3. Here, $B$ is a new independent variable.
\end{definition}

\begin{figure}
\centering 
\psfig{file=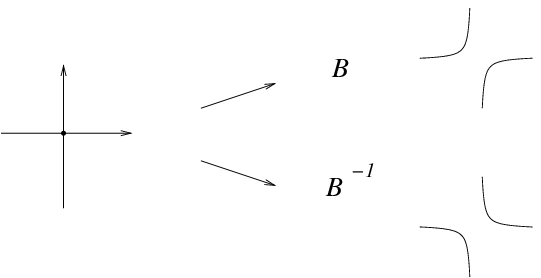}
\caption{}
\end{figure}

At crossings we smooth the diagram as usual in the definition of the Kauffman bracket. The {\em singular Kauffman 
bracket}, denoted $<D>_s$,
is then defined exactly as the usual Kauffman bracket besides the changing in the variables at the double points.
Let $w(D)$ be the writhe of $D$. 

\begin{definition}
$V^s_L(A,B) = (-A)^{-3w(D)}<D>_s \in \mathbb{Z}[A ,A^{-1},B]$. Here,  L denotes the singular link represented by D.

\end{definition}

Let $T$ be a Kauffman state for the Alexander polynomial (compare \cite{F}).

\begin{definition}
The contribution of a dot in the adjacent regions to a double point $p$ are shown in Fig. 4. Here, $B$ is a new independent 
variable.
\end{definition}

\begin{figure}
\centering 
\psfig{file=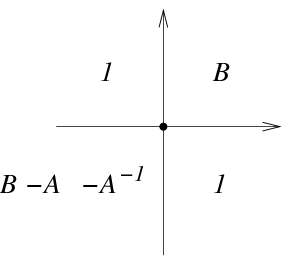}
\caption{}
\end{figure}

\begin{figure}
\centering 
\psfig{file=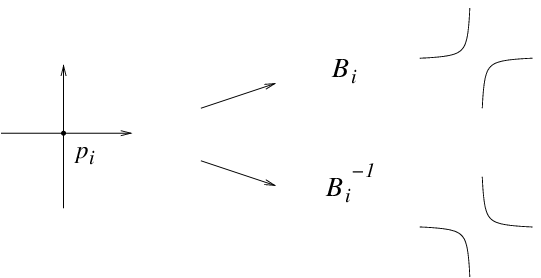}
\caption{}
\end{figure}

\begin{figure}
\centering 
\psfig{file=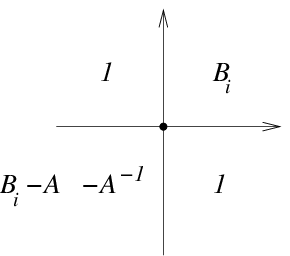}
\caption{}
\end{figure}

\begin{figure}
\centering 
\psfig{file=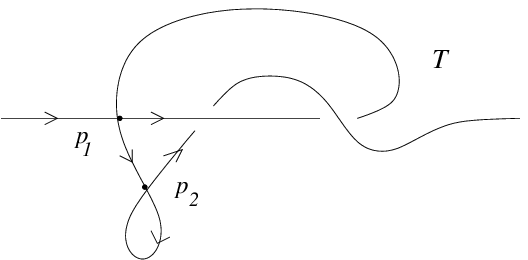}
\caption{}
\end{figure}

\begin{figure}
\centering 
\psfig{file=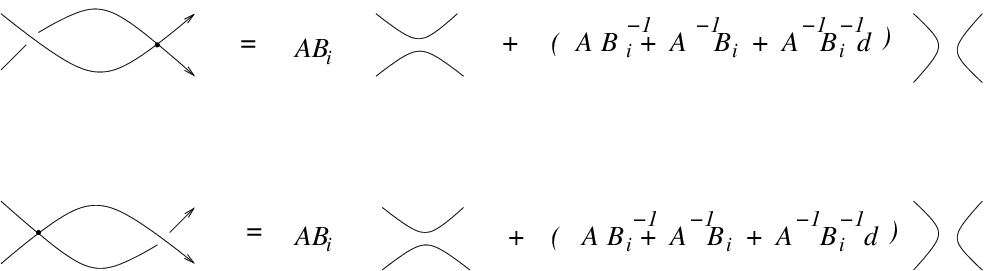}
\caption{}
\end{figure}

\begin{figure}
\centering 
\psfig{file=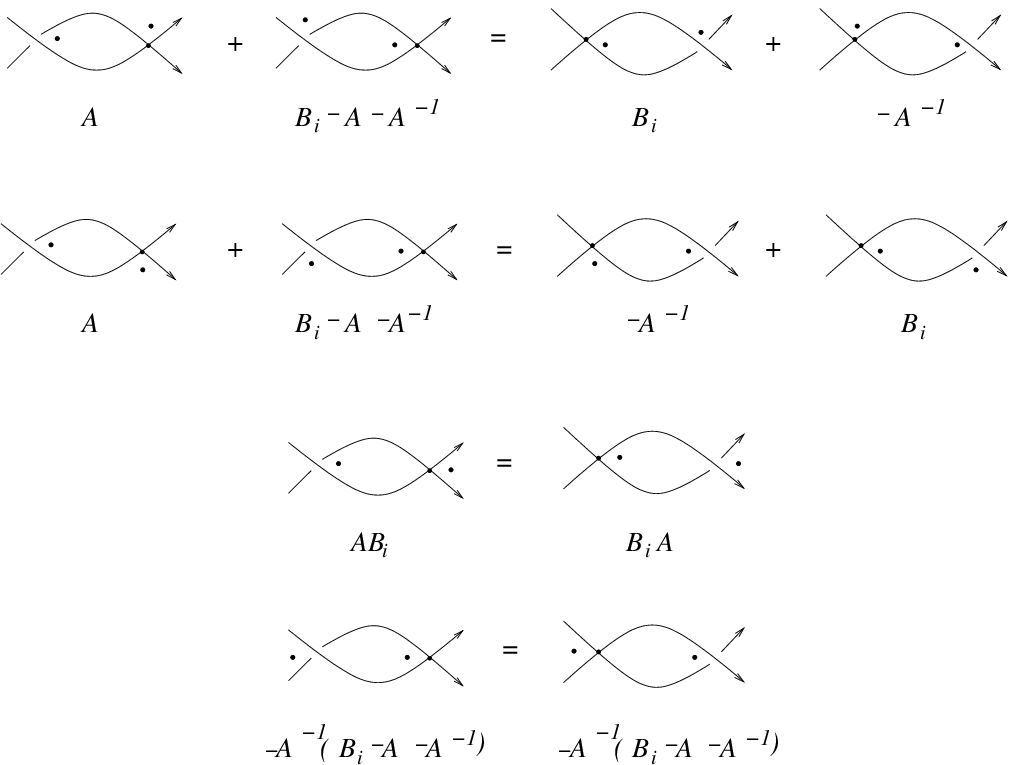}
\caption{}
\end{figure}

\begin{figure}
\centering 
\psfig{file=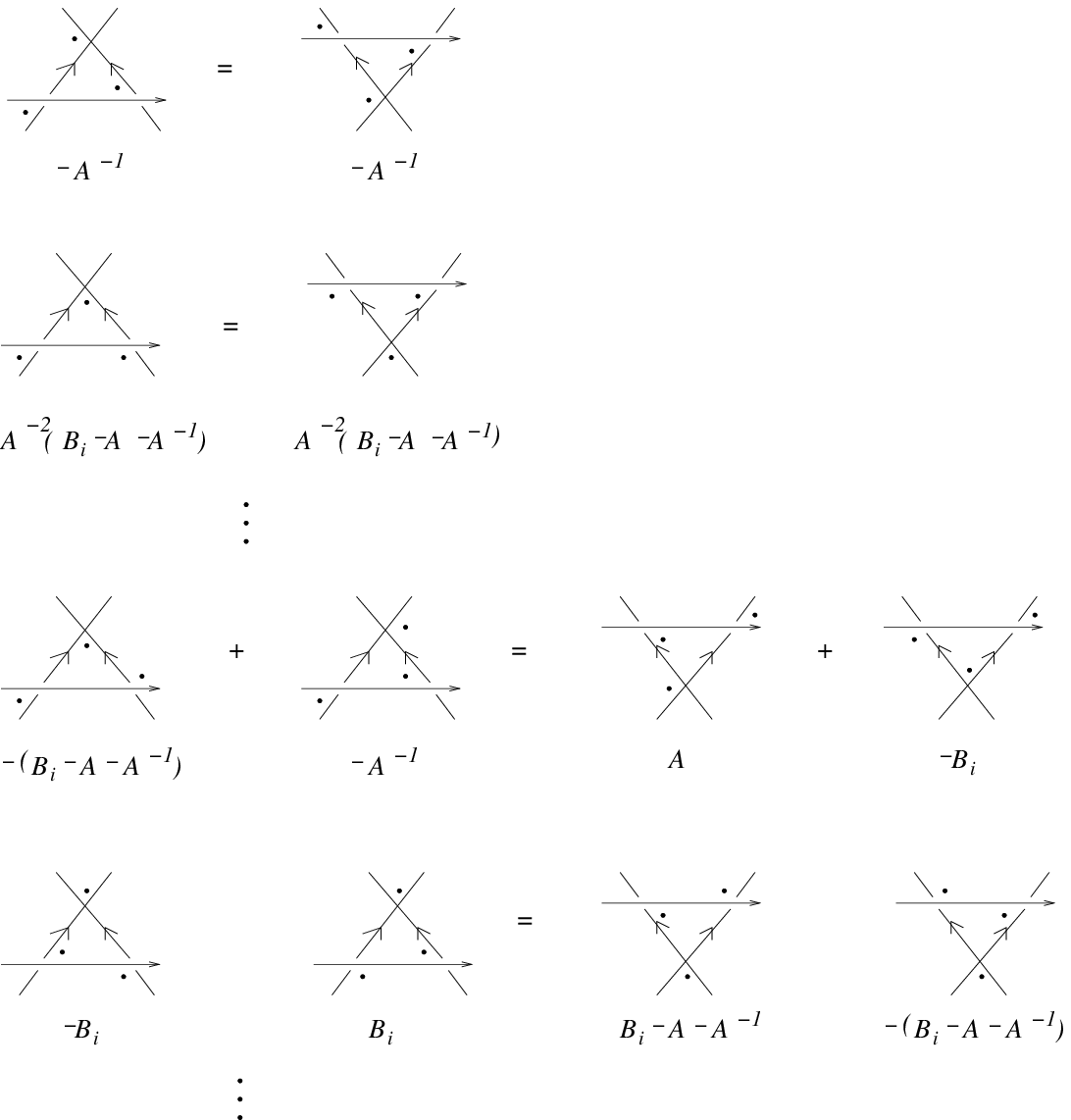}
\caption{}
\end{figure}

\begin{figure}
\centering 
\psfig{file=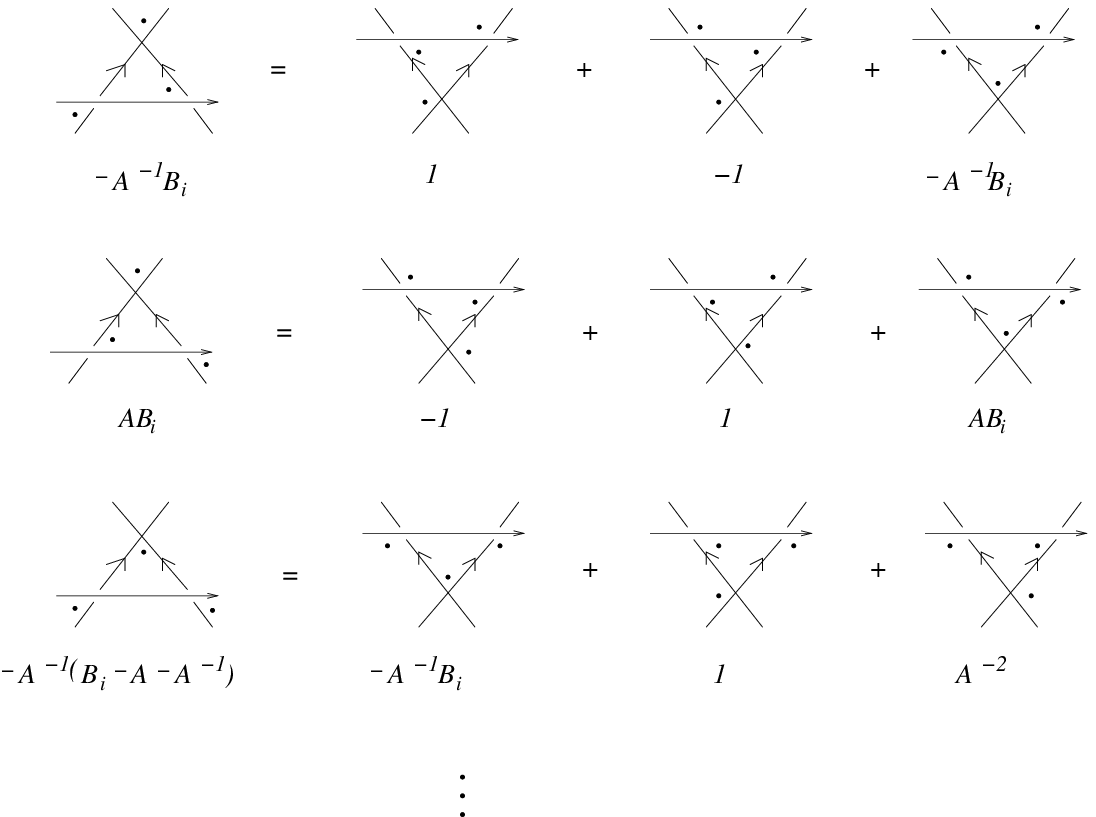}
\caption{}
\end{figure}

Notice, that the middle of Fig. 6 in \cite{OSS} coincides with our Fig. 4  for $B=0$.

To each state T we associate now the product of the monomials corresponding to the dots. 

\begin{definition}
The {\em Alexander 
polynomial for singular links\/}
 $\Delta^s_L(A,B) \in \mathbb{Z}[A ,A^{-1},B]$ is then the sum of all products of monomials over all possible states T.
 Here,  L denotes the singular link represented by D.
\end{definition}

\begin{theorem}
$V^s_L(A,B)$ and $\Delta^s_L(A,B)$ are isotopy invariants of oriented singular links (with a distinguished component
in the case of $\Delta^s_L(A,B)$).
\end{theorem}

\begin{remark}
Our state models for $V^s_L(A,B)$ and $\Delta^s_L(A,B)$ can be generalized in the evident way to define invariants for
singular links in handle bodies. We consider the handle body as a trivial line bundle over a planar surface. In order to define
the generalization of $\Delta^s_L(A,B)$ we have just to choose two boundary components of the planar surface and we mark
the adjacent regions of the complement of the link projection with the two stars (compare \cite{F}).
\end{remark}

Let us consider now long singular knots $K$ oriented from the left to the right (for their definition see e.g. \cite{V}). 
We assume that the knot $K$ has exactly d double point singularities.
Going along the
 knot there is a well defined order of the double point singularities induced by the first passage through each double 
point. We make the following trivial observation: this order is preserved under all isotopies of the long singular knot.
Indeed, each of the elementary moves involves at most one double point and can not change this order.

Let $p_i$ be the i-th double point. We replace the Figures 3 and 4 by the Figure 5 and respectively Figure 6 now.
Here $B_i$ is a new independent variable. We define the state sums now as previously but using the new variables
$B_i$ instead of the variable $B$. The result are two polynomials $V^s_K = (-A)^{-3w(D)}<D>_s
 \in \mathbb{Z}[A ,A^{-1},B_1,B^{-1}_1...B_d,B^{-1}_d]$ and 
$\Delta^s_K \in \mathbb{Z}[A ,A^{-1},B_1,...B_d]$. The polynomial $V^s_K$ is of degree 1 or -1 for each of the variables 
$B_i$ and $\Delta^s_K$ is of degree 0 or 1 for each of the variables $B_i$.  One could consider each of them as an ordered set of $2^d$ Laurent polynomials of the variable $A$.

\begin{theorem} The polynomials $V^s_K$ and $\Delta^s_K$ are isotopy invariants for oriented singular long knots. 
\end{theorem} 

\begin{remark}
If we consider knots instead of long knots then the variables $B_1,...B_d$ are only well defined up to 
permutations. Hence, we obtain in general only a partially ordered set of $2^d$ polynomials.
\end{remark}

\begin{example}
Let $T$ be the long singular trefoil as shown in Figure 7. One easily calculates 

$V^s_T = A^{-6}B^{-1}_1B^{-1}_2 + (-A^4-A^{-4})B_1B^{-1}_2 +(-A^{-4}-A^{-8})B^{-1}_1B_2 + (A^6+A^2+A^{-2}+A^{-6})B_1B_2$

$\Delta^s_T = A^{-2} + (A-A^{-1})B_1$.

The {\em inverse long knot} $-K$ is defined by reflecting the diagram of $K$ at a vertical line (which interchanges the 
endpoints) and then reversing the orientation on the knot. If we consider a long knot as a knot in the 3-sphere
then the inverse long knot corresponds just to the usual inverse knot in the 3-sphere. A long singular knot $K$ 
is called {\em invertible} if $K$ is isotopic to $-K$ as long singular knots.

The calculation of the polynomial $V^s_K$  for $-T$ is identical with that for $T$ besides the fact that the variables $B_1$
and $B_2$ have to be interchanged. Consequently, if $T$ would be  invertible then $V^s_T(A,B_1,B_2) = V^s_T(A,B_2,B_1)$.
This is not the case.

One easily calculates that

$\Delta^s_{-T} = A^{-2} + (A-A^{-1})B_2$.

Consequently, each of the two new polynomials shows that $T$ is not invertible. Notice, that the previous polynomials
  $V^s_L(A,B)$ and $\Delta^s_L(A,B)$ even applied to long 
knots fail to do this (as well as the HOMFLY-PT polynomial for singular links \cite{KV}).

\end{example}

\section{Proofs}

We indicate only the changings with respect to the proofs of the invariance for $V_L(A)$ and $\Delta_L(A)$.

We verify the invariance of $V^s_K$ under one of the moves $SII$ in Fig. 8 (all other cases are analogous). 

The proof that it is invariant under the 
moves $SIII$ is word for word identical with Kauffman's original proof (see \cite{K2}) that $V_L(A)$ is invariant under 
Reidemeister moves of type $RIII$. We have just to replace the crossing which Kauffman smoothes by the double 
point $p_i$.

We verify invariance of $\Delta^s_K$ under one of the moves $SII$ in Fig. 9 (all other cases are analogous).

To prove invariance under the moves $SIII$ we have to consider all Kauffman states $T$ near a triple point 
(compare Fig. 15 in \cite{F}). We carry out just one example in Fig. 10 and Fig. 11. The verification in all other cases is analogous 
and is left to the reader.

Theorem 1 is obtained by identifying all the variables $B_i$.

The theorems are proven.

Laboratoire de Math\'ematiques

Emile Picard

Universit\'e Paul Sabatier

118 ,route de Narbonne 

31062 Toulouse Cedex 09, France

fiedler@picard.ups-tlse.fr

\end{document}